\documentclass[12pt]{article}
\usepackage{amsmath}
\usepackage{amsbsy}
\usepackage{amsfonts}
\usepackage{mathtools}
\usepackage{bbm}
\usepackage{amsmath}
\usepackage{amstext}
\usepackage{amsthm}
\usepackage{bbold}
\usepackage{amssymb}
\usepackage[dvips]{graphicx}
\usepackage{amsfonts}
\usepackage{caption2}
\usepackage{amsfonts,mathrsfs}
\usepackage{makeidx}
\usepackage{fancyhdr}
\usepackage{epic,eepic,epsfig}
\usepackage{mdwlist}
\usepackage{amscd}
\usepackage{amsbsy}
\usepackage{epsfig,graphicx}
\usepackage{epsfig,mathrsfs,wasysym,stmaryrd}
\usepackage[all]{xypic}
\usepackage[knot,poly]{xy}
\usepackage{latexsym}
\usepackage{comment}

\theoremstyle{plain}

\newtheorem{theorem}{Theorem}

\newtheorem{lemma}{Lemma}

\newtheorem{corollary}{Corollary}
\newenvironment{pf}{\medskip\noindent{Proof:}
  \hspace{-.5cm}      \enspace}{\hfill \qed \newline \smallskip}

\date{}

\begin{document}

\begin{center}
\textbf{\LARGE{ The number of spanning trees in circulant graphs, its  arithmetic properties and asymptotic }}
\vspace{12pt}

{\large\textbf{A.~D.~Mednykh,}}\footnote{{\small\em Sobolev Institute of Mathematics,
Novosibirsk State University, smedn@mail.ru}}
{\large\textbf{I.~A.~Mednykh,}}\footnote{{\small\em Sobolev Institute of Mathematics,
Novosibirsk State University, ilyamednykh@mail.ru}}
\end{center}

\section*{Abstract}

In this paper, we develop a new method to produce explicit formulas for the number
$\tau(n)$ of spanning trees in the undirected circulant graphs $C_{n}(s_1,s_2,\ldots,s_k)$
and $C_{2n}(s_1,s_2,\ldots,s_k,n).$
Also, we prove that in both cases the number of spanning trees can be represented in the form
$\tau(n)=p \,n \,a(n)^2,$ where $a(n)$ is an integer sequence and $p$ is a prescribed natural
number depending   on the parity of $n.$ Finally, we find an asymptotic formula for $\tau(n)$
through the Mahler measure of the associated Laurent polynomial
$L(z)=2k-\sum\limits_{i=1}^k(z^{s_i}+z^{-s_i}).$
\bigskip

\noindent
\textbf{Key Words:} spanning tree, circulant graph, Laplacian matrix, Chebyshev polynomial,
Mahler measure\\\textbf{AMS classification:} 05C30, 39A10\\

\section{Introduction}

The \textit{complexity} of a finite connected graph $G$, denoted by $\tau(G),$ is the number
of spanning trees of $G.$ One of the first results on the complexity was obtained by Cayley \cite{Cay89}
who proved that the number of spanning trees in the complete graph $K_n$ on $n$ vertices is $n^{n-2}.$

The famous Kirchhoff's Matrix Tree Theorem~\cite{Kir47} states that $\tau(G)$ can be expressed
as the product of nonzero Laplacian eigenvalues of $G$ divided by the number of its vertices.
Since then, a lot of papers devoted to the complexity of various classes of graphs were published.
In particular, explicit formulae were derived for complete multipartite graphs~\cite{Cay89, Austin60},
almost complete graphs~\cite{Wein58}, wheels~\cite{BoePro}, fans~\cite{Hil74}, prisms~\cite{BB87},
ladders~\cite{Sed69}, M\"obius ladders~\cite{Sed70}, lattices \cite{SW00}, anti-prisms \cite{SWZ16},
complete prisms~\cite{Sch74} and for many other families. For the circulant graphs some explicit and recursive
formulae are given in \cite{Xiebin, XiebinLinZhang, YTA97, ZY99, ZhangYongGol, ZhangYongGolin}.

Starting with Boesch and Prodinger \cite{BoePro} the idea to study the complexity of graphs by making
use of Chebyshev polynomials was implemented. This idea provided a way to find complexity of circulant
graphs and their natural generalisations in~\cite{KwonMedMed, Med1, Xiebin, XiebinLinZhang, YTA97,   
ZhangYongGol, ZhangYongGolin}.

Recently, asymptotical behavior of complexity for some families of graphs was investigated from 
the point of view of so called Malher measure \cite{GutRog}, \cite{SilWil}, \cite{SilWil1}. 
Mahler measure of a polynomial $P(z)$, with complex coefficients, is the product of the roots 
of $P(z)$ whose modulus is greater than $1$ multiplied by the leading coefficient. For general 
properties of the Mahler measure see survey \cite{Smyth08} and monograph \cite{EverWard}. 
It worth mentioning that the Mahler measure is related to the growth of groups, values of some 
hypergeometric functions and volumes of hyperbolic manifolds \cite{Boy02}.

For a sequence of graphs $G_n$ with the number of vertices $v(G_n),$ one can consider
the number of spanning trees $\tau(G_n)$ as a function of $n.$ Assuming that the limit
$\lim_{n\to \infty}\frac{\log\tau(G_n)}{v (G_n)}$ exists, it is sometimes called the
associated tree entropy or the thermodynamic limit of the family $G_n$ \cite{Lyon05}.
This number plays an important role in statistical physics and was investigated by many
authors (\cite{Kas2}, \cite{SW00}, \cite{SilWil}, \cite{SilWil1}, \cite{TF}, \cite{Wu77}).

The purpose of this paper is to present new formulas for the number of spanning trees
in circulant graphs and investigate their arithmetical properties and asymptotic. We
mention that the number of spanning trees for circulant graphs was found earlier in
(\cite{Louis}, \cite{Xiebin}, \cite{XiebinLinZhang}, \cite{ZhangYongGol}, \cite{ZhangYongGolin})
.

The structure of the paper is as follows. First, in the sections~\ref{count} and \ref{oddcomplexity}
we present new explicit formulas for the number of spanning trees in the undirected
circulant graphs $C_{n}(s_1,s_2,\ldots,s_k)$ and $C_{2n}(s_1,s_2,\ldots,s_k,n)$
of even and odd valency respectively. They will be given in terms of Chebyshev polynomials. These formulas are
different from those obtained earlier in the papers (\cite{Xiebin}, \cite{XiebinLinZhang}, \cite{ZhangYongGol},
\cite{ZhangYongGolin}). Moreover, by our opinion, the obtained
formulas are more convenient for analytical investigation. Next, in the section~\ref{circarithm}
we provide some arithmetic properties of the complexity function. More precisely, we show
that the number of spanning trees of the  circulant graph  can be represented in
the form $\tau(n)=p \,n \,a(n)^2,$ where $a(n)$ is an integer sequence and $p$ is a prescribed
natural number depending   on the parity of $n.$ Later, in the section~\ref{assection}, we use
explicit formulas for the complexity in order to produce its asymptotic in terms of Mahler
measure of the associated polynomials. For circulant graphs of even valency   the associated polynomial is
$L(z)=2k-\sum\limits_{j=1}^{k}(z^{s_j}+z^{-s_j}).$ In this case  (Theorem~\ref{asymptotic1}),  assuming
${\rm gcd}\,(s_1,s_2,\ldots,s_k)=d,$  we have $\tau(n)\sim\frac{ n\, d^2}{q}A^n,\,n\to\infty,$  where
$q=s_1^2+s_2^2+\ldots+s_k^2,\,A=M(L)$  and  $M(L)$ stands for the Mahler measure of $L(z).$
For circulant graphs of odd valency
we use the polynomial $R(z)=L(z)(L(z)+2).$
Then the respective asymptotic (Theorem~\ref{asymptotic3}) is $\tau(n)\sim\frac{n \,d^2}{2\,q}K^n,\,n\to\infty,$  where
$K=M(R).$   As a consequence (Corollary~\ref{cor1} and Corollary~\ref{cor2}), we obtained that the thermodynamic limits
of the sequences $C_{n}(s_1,s_2,\ldots,s_k)$ and $C_{2n}(s_1,s_2,\ldots,s_k,n)$ are
$\log M(L)$ and $\log M(R)$ respectively. In the last section~\ref{tables},
we illustrate the obtained results by a series of examples.

\section{Basic definitions and preliminary facts}

Consider a finite connected  graph $G$  without
loops. We denote the vertex and edge set of $G$ by $V(G)$ and $E(G),$ respectively.
Given $u, v\in V(G),$ we set $a_{uv}$ to be equal to the number of edges between
vertices $u$ and $v.$ The matrix $A=A(G)=\{a_{uv}\}_{u, v\in V(G)}$ is called
\textit{the adjacency matrix} of the graph $G.$ The degree $d(v)$ of a vertex
$v \in V(G)$ is defined by $d(v)=\sum_{u\in V(G)}a_{uv}.$ Let $D=D(G)$ be the
diagonal matrix indexed by the elements of $V(G)$ with $d_{vv} = d(v).$ The matrix
$L=L(G)=D(G)-A(G)$ is called \textit{the Laplacian matrix}, or simply \textit{Laplacian},
of the graph $G.$  

By $I_n$ we denote the identity matrix of order $n.$

Let $s_1,s_2,\ldots,s_k$ be integers such that $1\leq s_1<s_2<\ldots<s_k\leq\frac{n}{2}.$
The graph $C_{n}(s_1,s_2,\ldots,s_k)$ with $n$ vertices $0,1,2,\ldots,~{n-1}$ is called
\textit{circulant graph} if the vertex $i,\, 0\leq i\leq n-1$ is adjacent to the vertices
$i\pm s_1,i\pm s_2,\ldots,i\pm s_k\ (\textrm{mod}\ n).$ When $s_k<\frac{n}{2}$ all vertices of
a graph have even degree $2k.$ If $n$ is even and $s_k=\frac{n}{2},$ then all vertices have odd degree $2k-1.$
It is well known  that the circulant $C_n(s_1,s_2,\ldots,s_k)$ is connected if and
only if ${\rm gcd}\,(s_1,s_2,\ldots,s_k,n)~=~1. $  More generally, the number of connected
components of $C_n(s_1,s_2,\ldots,s_k)$ is $d=\textrm{gcd}\,(s_1,s_2,\ldots,s_k,n),$
with each of the vertices $0,1,...,d - 1$ lying in different components, and with each
component being isomorphic to $C_{n/d}(s_1/d,s_2/d,\ldots,s_k/d).$ So, for $d>1$ graph
is disconnected and has no spanning trees. In what follows, all graphs are supposed to
be connected.

Two circulant graphs $C_{n}(s_1,s_2,\ldots,s_k)$ and $C_{n}(\tilde{s}_1,\tilde{s}_2,\ldots,\tilde{s}_k)$
of the same order are said to be conjugate by multiplier if there exists an integer $r$
coprime to $n$ such that $\{\tilde{s}_1,\tilde{s}_2,...,\tilde{s}_k\}=\{rs_1,rs_2,\ldots,rs_k\}$
as subsets of $\mathbb{Z}_n.$ In this case, the graphs are isomorphic, with multiplication
by the unit $ r\,(\textrm{mod}\,n)$  giving an isomorphism.

In 1967, A.~\'Ad\'am conjectured that two circulant graphs are isomorphic if and only if they are
conjugate by a multiplier \cite{Adam}. It should be mentioned that his original goal was to give
a classification of all isomorphism classes of circulant graphs. However, there is the following
counterexample to \'Ad\'amÕs Conjecture. The graphs $C_{16}(1,2,7)$ and $C_{16}(2,3,5)$ are isomorphic,
but they are not conjugate by a multiplier \cite{CondGra}. A complete solution of the isomorphism
problem for circulant graphs was done by M.~Muzychuk  \cite{Muz}.

We call an $n\times n$ matrix \textit{circulant,} and denote it by $circ(a_0, a_1,\ldots,a_{n-1})$
if it is of the form
$$circ(a_0, a_1,\ldots, a_{n-1})=
\left(\begin{array}{ccccc}
a_0 & a_1 & a_2 & \ldots & a_{n-1} \\
a_{n-1} & a_0 & a_1 & \ldots & a_{n-2} \\
  & \vdots &   & \ddots & \vdots \\
a_1 & a_2 & a_3 & \ldots & a_0\\
\end{array}\right).$$

It easy to see that adjacency and Laplacian matrices of  the circulant graph are circulant
matrices. The converse is also true. If the Laplacian matrix of a graph is circulant then
the graph is also circulant.

Recall \cite{PJDav} that the eigenvalues of matrix $C=circ(a_0,a_1,\ldots,a_{n-1})$ are
given by the following simple formulas $\lambda_j=L(\varepsilon^j_n),\,j=0,1,\ldots,n-1,$
where $L(x)=a_0+a_1 x+\ldots+a_{n-1}x^{n-1}$ and $\varepsilon_n$ is an order $n$ primitive
root of the unity. Moreover, the circulant matrix $C=L(T),$ where $T=circ(0,1,0,\ldots,0)$
is the matrix representation of the shift operator
$T:(x_0,x_1,\ldots,x_{n-2},x_{n-1})\rightarrow(x_1, x_2,\ldots,x_{n-1},x_0).$  

Let $P(z) = a_0 z^d +\ldots+a_d = a_0 \prod\limits_{i=1}^d(z-\alpha_i)$ be a nonconstant polynomial
with  complex coefficients. Then, following Mahler \cite{Mahl62} its {\it Mahler measure} is defined to be
\begin{equation} M(P):=\exp(\int_0^1\log|P(e^{2\pi i t})|dt),\end{equation}
the geometric mean of $|P (z)|$ for $z$ on the unit circle. However, $M (P)$ had appeared earlier in a paper by
Lehmer \cite{Lehm33}, in an alternative form

\begin{equation}  M(P)=\alpha_0\prod\limits_{|\alpha_i|>1}|\alpha_i|.\end{equation}
The equivalence of the two definitions follows immediately from Jensen�s formula \cite{Jen99}
\begin{equation*} \int_0^1\log|e^{2\pi i t}-\alpha|dt=\log_+|\alpha|,\end{equation*}
where $\log_+x$  denotes $\max(0,\,\log x).$ Sometimes, it is more convenient to deal with the \textit{small Mahler measure}
which is defined as $$m(P):=\log M(P)=\int_0^1\log|P(e^{2\pi i t})|dt.$$
The concept of Mahler measure can be naturally extended to the class of Laurent polynomials
$P(z) = a_{0}z^{p}+a_{1}z^{p+1}+\ldots+a_{s-1}z^{p+s-1}+a_{s}z^{p+s}= a_{s}z^p\prod\limits_{i=1}^{s}(z-\alpha_i),$ where $a_s\ne0$ and $p$ is an arbitrary integer (not necessarily positive).
The following  properties of Mahler measure can be easy deduced from definition. 
Given polynomials $P(z)$ and $Q(z)$ and a positive integer $d,$ we have 
$M(P(z)\,Q(z))=M(P(z)) M(Q(z))$ and $M(P(z^d))=M(P(z)).$

During the paper, we will use the  basic properties of Chebyshev polynomials. Let $T_n(z)=\cos(n\arccos z)$  and
$U_{n-1}(z)={\sin(n\arccos z)}/{\sin(\arccos z)}$ be the Chebyshev polynomials of the first and   second kind respectively.

Then $  T^{\prime}_n(z)=n\,U_{n-1}(z),\,T_n(1)=1,\,U_{n-1}(1)=n.$ For  $z\neq 0$  we have $T_n(\frac12(z+z^{-1}))=\frac12(z^n+z^{-n}),$  and for  $z\neq 1$
the identity   $(T_n(z)-1)/(z-1)=U_{n-1}^2(\sqrt{(z+1)/2})$  holds.

Also, $T_n(z)$ and $U_{n-1}(z)$ admit the following quantum representation
$T_n(z)=({q^n+q^{-n}})/{2}$ and $U_{n-1}(z)=({q^n-q^{-n}})/({q-q^{-1}}),$
where $q=z+\sqrt{z^2-1}.$ See \cite{MasHand} for more advanced properties.

\section{Complexity of circulant graphs of even valency}\label{count}

The aim of this section is to find  new formulas for the numbers of spanning trees of
circulant graph $C_{n}(s_1,s_2,\ldots,s_k)$ in terms of Chebyshev polynomials. It should
be noted that nearby results were obtained earlier by different methods in the papers
\cite{Xiebin}, \cite{ZhangYongGol}, \cite{XiebinLinZhang}, \cite{ZhangYongGolin}.

\begin{theorem}\label{theorem1}
The number of spanning trees $\tau(n)$ in the circulant graph $C_{n}(s_1,s_2,\ldots,s_k),$
$1\le s_1< s_2<\ldots< s_k<\frac{n}{2},$ is given by the formula
$$\tau(n)=(-1)^{(n-1)(s_k-1)}n\prod_{p=1}^{s_k-1}\frac{T_n(w_p)-1}{w_p-1},$$ thereby $w_p,\, p=1,2,\ldots,s_k-1$
are roots of the algebraic equation $P(w)=0,$ where
$$P(w)= \sum_{j=1}^k \frac{T_{s_j}(w)-1}{w-1}$$
and $T_k(w)$ is the Chebyshev polynomial of the first kind.
\end{theorem}

\medskip\noindent{Proof:} By the celebrated Kirchhoff theorem, the number of spanning trees $\tau(n)$ is
equal to the product of nonzero eigenvalues of the Laplacian of a graph $C_{n}(s_1,s_2,\ldots,s_k)$ divided
by the number of its vertices $n.$ To investigate the spectrum of Laplacian matrix, we denote by $T=circ(0,1,\ldots,0)$
the $n \times n$ shift operator. Consider the Laurent polynomial $L(z)=2k-\sum\limits_{i=1}^k(z^{s_i}+z^{-s_i}).$

Then the Laplacian of  $ C_{n}(s_1,s_2,\ldots,s_k)$ is given by the matrix
$$\mathbb{L}=L(T)=2k I_n-\sum\limits_{i=1}^k(T^{s_i}+T^{-s_i}).$$
The eigenvalues of circulant matrix $T$ are $\varepsilon_n^j,\,j=0,1,\ldots,n-1,$ where
$\varepsilon_n=e^\frac{2\pi i}{n}.$ Since all of them are distinct, the matrix $T$ is
conjugate to the diagonal matrix $\mathbb{T}=diag(1,\varepsilon_n,\ldots,\varepsilon_n^{n-1})$
with diagonal entries $1,\varepsilon_n,\ldots,\varepsilon_n^{n-1}$. To find spectrum of
$\mathbb{L},$ without loss of generality, one can assume that $T=\mathbb{T}.$
Then $\mathbb{L}$ is a diagonal matrix. This essentially simplifies the problem of
finding eigenvalues of $\mathbb{L}.$ Indeed, let $\lambda$ be an eigenvalue of $L$ and
$x$ be the respective eigenvector. Then we have the following system of linear equations
$$((2k-\lambda)I_n-\sum\limits_{i=1}^k(T^{s_i}+T^{-s_i}))x=0.$$

Recall the matrices under consideration are diagonal and the $(j+1,j+1)$-th entry of $T$ is
equal to $\varepsilon_n^{j},$ where $\varepsilon_n=e^{\frac{2\pi i}{n}}.$

Let ${\bf e}_{j} =(0,\ldots,\underbrace{1}_{j-th},\ldots, 0),\,j=1,\ldots,n.$
Then,  for any $j=0,\ldots, n-1,$ matrix $\mathbb{L}$ has an eigenvalue
$\lambda_j=L(\varepsilon_n^j)=2k-\sum\limits_{i=1}^k(\varepsilon_n^{j s_i}+\varepsilon_n^{-j s_i})$
with eigenvector ${\bf e}_{j+1}.$ Since the graph under consideration is connected,
we have $\lambda_0=0$ and $\lambda_j>0,\, j=1,2,\ldots, n-1.$ Hence
$$\tau(n)=\frac{1}{n}\prod\limits_{j=1}^{n-1}L(\varepsilon_n^j).$$
To continue the calculation of $\tau(n)$ we need the following  lemmas.

\begin{lemma}\label{lemma2} The following identity holds
$$L(z)=2(1-w) P(w),$$ where  $P(w)$ is an  integer polynomial given by the formula
$$P(w)= \sum_{i=1}^k \frac{T_{s_j}(w)-1}{w-1},$$
$T_k(w)$ is the Chebyshev polynomial of the first kind and $w=\frac{1}{2}(z+z^{-1}).$
\end{lemma}

\begin{pf}  Let us substitute $z=e^{i\varphi}.$ It is easy to
see that $w=\frac{1}{2}(z+z^{-1})=\cos\varphi,$ so we have
$T_k(w)=\cos(k\arccos w)=\cos(k\varphi).$ Also, $L(z)=2k-\sum_{i=1}^k \cos(k\varphi).$
Then the statement of the lemma follows from elementary calculations.
\end{pf}

   To find find asymptotic for the number of spanning trees we have to use the following property of roots of polynomial $L(z).$
See Lemma 6 in \cite{ZhangYongGol} for a similar result.

\begin{lemma}\label{lemma3} Suppose that ${\rm gcd}\,(s_1,s_2,\ldots,s_k)=1.$ Then the roots of Laurent  polynomial  
$L(z)=2k-\sum\limits_{i=1}^k(z^{s_i}+z^{-s_i})$  are
$1,1,z_1,1/z_1,\ldots,z_{s_k-1},1/z_{s_k-1},$ where for all $s=1,\ldots,s_k-1,|z_s|\neq 1$  and $w_s=\frac{1}{2}(z_s+z_s^{-1})$
is a root of polynomial $P(w).$
\end{lemma}

\begin{pf}By Lemma \ref{lemma2}, $L(z)=2(1-w) P(w),$ where $w=\frac{1}{2}(z+z^{-1})$ and
$P(w)= \sum\limits_{i=1}^k \frac{T_{s_j}(w)-1}{w-1}$ is the polynomial of degree $s_k-1.$
Note that $2(1-w) = -\frac{(z-1)^2}{z}.$ Since
$P(1)=\sum\limits_{i=1}^k T_{s_j}^{\prime}(1)=\sum\limits_{i=1}^k {s_j}^2\neq 0,$
the Laurent polynomial $L(z)$ has the root $z=1$ with multiplicity two. Hence, the roots of $L(z)$ are
$1,1,z_1,1/z_1,\ldots,z_{s_k-1},1/z_{s_k-1},$ where for all $s=1,\ldots,s_k-1,z_s\neq 1;$
the respective roots of $P(w)$ are $1, w_s=\frac{1}{2}(z_s+z_s^{-1}),\,s=1,\ldots,s_k-1.$

Now we show  that condition ${\rm gcd}\,(s_1,s_2,\ldots,s_k)=1$ implies   $|z_s|\neq 1$
for any $s=1,\ldots,s_k-1.$ Indeed, we already have $z_s\neq 1.$ Suppose that $|z_s|=1.$
Then for some real number $\varphi$  we have $z_s=e^{i\,\varphi}.$ Since $z_s$ is a root
of $L(z),$  we obtain
$$L(e^{i\,\varphi})=2k-\sum\limits_{i=1}^k(e^{i\,s_i\varphi}+
e^{-i\,s_i\varphi})=2\sum\limits_{i=1}^k(1-\cos(s_i\varphi))=0.$$
Hence, $\cos(s_i\varphi)=1$ and $s_i\varphi=2\pi m_i$ for some integers $m_i,\,i=1,2,\ldots,k.$         
By virtue of ${\rm gcd}\,(s_1,s_2,\ldots,s_k)=1,$ one can find integers $x_i,\,i=1,2,\ldots,k$
such that $x_1s_1+x_2s_2+\cdots+x_ks_k=1.$ (See, for example, \cite{Apost}, p.~21). Then
$$z_s=e^{i\,\varphi}=e^{i(x_1s_1+\cdots+x_ks_k)\varphi}=e^{i(x_1s_1\varphi+\cdots+
x_k s_k\varphi)}=e^{i(2\pi x_1m_1 +\cdots+2\pi ix_k s_k)}=1.$$ Contradiction.
\end{pf}
\medskip

\begin{lemma}\label{lemma5}
Let $H(z)=\prod\limits_{s=1}^{m}(z-z_s)(z-z_s^{-1})$ and $H(1)\neq0.$ Then
$$\prod\limits_{j=1}^{n-1}H(\varepsilon_n^j)=\prod\limits_{s=1}^{m}\frac{T_n(w_s)-1}{w_s-1},$$
where $w_s=\frac12(z_s+z_s^{-1}),\,s=1,\ldots,k$ and $T_n(w)$ is the
Chebyshev polynomial of the first kind.
\end{lemma}

\begin{pf}
It is easy to check that $\prod\limits_{j=1}^{n-1}(z-\varepsilon_n^j)=\frac{z^n-1}{z-1}$
if $z\neq1.$ Also we note that $\frac12(z^n+z^{-n})=T_n(\frac12(z+z^{-1})).$ By the
substitution $z=e^{i\,\varphi}$ the latter follows from the evident identity
$\cos(n\varphi)=T_n(\cos\varphi).$ Then we have
\begin{eqnarray*}
\prod\limits_{j=1}^{n-1} H(\varepsilon_n^j) &=& \prod\limits_{j=1}^{n-1}\prod\limits_{s=1}^{m}(\varepsilon_n^j-z_s)(\varepsilon_n^j-z_s^{-1})\\
&=&\prod\limits_{s=1}^{m}\prod\limits_{j=1}^{n-1}(z_s-\varepsilon_n^j)(z_s^{-1}-\varepsilon_n^j)\\
&=&
\prod\limits_{s=1}^{m}\left(\frac{z_s^n-1}{z_s-1}\cdot\frac{z_s^{-n}-1}{z_s^{-1}-1}\right)=
\prod\limits_{s=1}^{m}\frac{T_n(w_s)-1}{w_s-1}.\vspace{-1cm}
\end{eqnarray*}
\end{pf}
\bigskip

To continue the proof of the theorem  we set $H(z)=\prod\limits_{s=1}^{s_k-1}(z-z_{s})(z-z_{s}^{-1}).$
Then $L(z)=-\frac{(1-z)^2}{z^{s_k}}H(z).$

Note that
$\prod\limits_{j=1}^{n-1}(1-\varepsilon_n^j)=\lim\limits_{z\to1}\prod\limits_{j=1}^{n-1}
(z-\varepsilon_n^j)= \lim\limits_{z\to1}\frac{z^n-1}{z-1}=n$ and
$\prod\limits_{j=1}^{n-1}\varepsilon_n^{j} = (-1)^{n-1}$.  As a
result, taking into account Lemma~\ref{lemma5}, we obtain
\begin{eqnarray*}
\tau(n)&=&\frac{1}{n}\prod\limits_{j=1}^{n-1}L(\varepsilon_n^j)=\frac{1}{n}\prod\limits_{j=1}^{n-1}
(-\frac{(1-\varepsilon_n^j)^2}{(\varepsilon_n^{j})^{s_k}}H(\varepsilon_n^j))=\frac{(-1)^{(n-1)(s_k-1)}n^2}{n}
\prod\limits_{j=1}^{n-1}H(\varepsilon_n^j)\\
&=&
(-1)^{(n-1)(s_k-1)}n\prod\limits_{s=1}^{s_k-1}\frac{T_n(w_s)-1}{w_s-1}.
\end{eqnarray*}
The theorem is proved.
$\hfill \qed$

\bigskip
The next corollary gives an important tool to find asymptotic behavior for the number of spanning trees.
It will be done later in  section~\ref{assection}.
\bigskip

\begin{corollary}\label{corollary21} The  number 
 of spanning trees in the circulant graph $C_{n}(s_1,s_2,\ldots,s_k) $  is given by the formula 
\begin{equation}\label{newform}  \tau(n)=\frac{n}{q}
\prod_{p=1}^{s_k-1}|2\,T_n(w_p)-2|,
\end{equation}
 where   $q=s_1^2+s_2^2+\ldots+s_k^2$  and $w_p, \,p=1,2,\ldots,s_k-1$  are different from $1$   roots of the equation  $\sum_{j=1}^k T_{s_j}(w)=k.$  
\end{corollary}
\bigskip

\begin{pf} We note that $w_p,\, p=1,2,\ldots,s_k-1$ are all roots of the polynomial
$Q(w)=(w-1)P(w)= \sum_{j=1}^k (T_{s_j}(w)-1)$ different from $1.$ We have
$Q^{\prime}(1)= s_1^2+s_2^2+\ldots+s_k^2=q.$ Since $Q(1)=0$ and $Q(w)$ is the order
$s_k$ polynomial whose leading term is $2^{s_k-1} w^{s_k},$ we have
$\prod_{p=1}^{s_k-1}(w_p-1)= (-2)^{1-s_k} Q^{\prime}(1).$ Taking into
account these properties, by Theorem~\ref{theorem1} we obtain
\begin{equation}\label{oldform}\tau(n)=\frac{(-1)^{(n-1)(s_k-1)}n}{(-2)^{1-s_k}Q^{\prime}(1)}\prod\limits_{s=1}^{s_k-1}
(T_n(w_s)-1)=\frac{(-1)^{n (s_k-1)} n\, 2^{s_k-1} }{q}\prod_{p=1}^{s_k-1}(T_n(w_p)-1).\end{equation}
Since $\tau(n)$ is a positive integer, (\ref{oldform}) implies (\ref{newform}).
\end{pf}

\begin{corollary}\label{corollary22}
$\tau(n)=n \left|\prod\limits_{s=1}^{s_k-1}U_{n-1}(\sqrt{\frac{w_p+1}{2}})\right|^2,$
where $w_p,\,p=1,2,\ldots,s$ are the same as above and $U_{n-1}(w)$
is the Chebyshev polynomial of the second kind.
\end{corollary}

\begin{pf} Follows from the identity $\frac{T_n(w)-1}{w-1}=U_{n-1}^2(\sqrt{\frac{w +1}{2}}).$ \end{pf}

\section{Complexity of circulant graphs of odd valency}\label{oddcomplexity}

The aim of this section is to find a new formula for the numbers of spanning trees
of circulant graph $C_{2n}(s_1,s_2,\ldots,s_k,n)$ in terms of Chebyshev polynomials.
Notice  that nearby results were obtained earlier by different methods in the papers \cite{Louis}, \cite{XiebinLinZhang},
\cite{ZhangYongGol},  \cite{ZhangYongGolin}.

\begin{theorem}\label{odddegree}
Let $C_{2n}(s_1,s_2,\ldots,s_k,n),\,1\leq s_1<s_2<\ldots<s_k<n,$ be a circulant graph of odd
degree. Then the number $\tau(n)$ of spanning trees in the graph $C_{2n}(s_1,s_2,\ldots,s_k,n)$
is given by the formula
$$\tau(n)= \frac{n\,4^{s_k-1}}{q}\prod_{p=1}^{s_k-1}(T_n(u_p)-1)\prod_{p=1}^{s_k}(T_n(v_p)+1),$$
where $q=s_1^2+s_2^2+\ldots+s_k^2,$ the numbers $u_p,p=1,2,\ldots,s_k-1 \text{ and } v_p,\,p=1,2,\ldots,s_k$
are respectively the roots of the algebraic equations $P(u)-1=0,\,u\ne1$ and $P(v)+1=0,$
where $P(w)=2k+1-2\sum\limits_{i=1}^{k}T_{s_i}(w)$ and $T_k(w)$ is the Chebyshev polynomial
of the first kind.
\end{theorem}

\begin{pf}
The Laplace operator of the graph $C_{2n}(s_1,s_2,\ldots,s_k,n)$ can be represented in the form
$$\mathbb{L}=(2k+1)I_{2n}-\sum\limits_{j=1}^{k}(T^{s_j}+T^{-s_j})-T^n,$$ where $T$ is $2n\times 2n$
shift  operator satisfying the equality $T^{2n}=I_{2n}.$ The eigenvalues of circulant matrix $T$
are $\varepsilon_{2n}^j,\,j=0,1,\ldots,2n-1,$ where $\varepsilon_{2n}=e^\frac{2\pi i}{2n}.$ Since
all of them are distinct, the matrix $T$ is conjugate to the diagonal matrix
$\mathbb{T}=diag(1,\varepsilon_{2n},\ldots,\varepsilon_{2n}^{2n-1})$ with diagonal entries
$1,\varepsilon_{2n},\ldots,\varepsilon_{2n}^{2n-1}$. To find spectrum of $\mathbb{L},$ without
loss of generality, one can assume that $T=\mathbb{T}.$ Then
$\mathbb{L}=diag(\lambda_0,\lambda_1,\ldots,\lambda_{2n-1})$ is the diagonal matrix with eigenvalues
$$\lambda_j=2k+1-\sum\limits_{l=1}^{k}(\varepsilon_{2n}^{j\, s_l}+
\varepsilon_{2n}^{-j \,s_l})-\varepsilon_{2n}^{j n}, \,j=0,1,\ldots, 2n-1.$$
Consider the following Laurent polynomial $L(z)=2k+1-\sum\limits_{i=1}^{k}(z^{s_i}+z^{-s_i}).$
Since $\varepsilon_{2n}^n=-1,$  we can write $\lambda_j=L(\varepsilon_{2n}^{j})-1$ if $j$ is even
and $\lambda_j=L(\varepsilon_{2n}^{j})+1$ if $j$ is odd. We note that $\lambda_0=0,$ so all
non-zero eigenvalues of $\mathbb{L}$  are $\lambda_1,\lambda_2,\ldots,\lambda_{2n-1}.$
By the Kirchhoff theorem we have
$$\tau(n)=\frac{\lambda_1\lambda_2\cdots\lambda_{2n-1}}{2n}=\frac{1}{2n}\prod\limits_{s=1}^{n-1}
(L(\varepsilon_{2n}^{2s})-1)\prod\limits_{s=0}^{n-1}(L(\varepsilon_{2n}^{2s+1})+1)$$

$$=\frac{1}{2n}\prod\limits_{s=1}^{n-1}(L(\varepsilon_{2n}^{2s})-1)\frac{\prod\limits_{p=1}^{2n-1}
(L(\varepsilon_{2n}^{p})+1)}{\prod\limits_{s=1}^{n-1}(L(\varepsilon_{2n}^{2s})+1)}=
\frac{1}{2n}\prod\limits_{s=1}^{n-1}(L(\varepsilon_{n}^{s})-1)\frac{\prod\limits_{p=1}^{2n-1}
(L(\varepsilon_{2n}^{p})+1)}{\prod\limits_{s=1}^{n-1}(L(\varepsilon_{n}^{s})+1)}.$$
\end{pf}

By making use of Lemma \ref{lemma5} and arguments from the proof of Theorem \ref{theorem1} we obtain \begin{enumerate}
\item[(i)]
$\prod\limits_{s=1}^{n-1}(L(\varepsilon_{n}^{s})-1)=(-1)^{(n-1)(s_k-1)}n^2
\prod_{p=1}^{s_k-1}\frac{T_n(u_p)-1}{u_p-1},$

\item[(ii)]
$\prod\limits_{s=1}^{n-1}(L(\varepsilon_{n}^{s})+1)=(-1)^{(n-1)
(s_k-1)}\prod_{p=1}^{s_k}\frac{T_n(v_p)-1}{v_p-1},$  and

\item[(iii)]
$\prod\limits_{p=1}^{2n-1}(L(\varepsilon_{2n}^{p})+1)=(-1)^{(2n-1)
(s_k-1)}\prod_{p=1}^{s_k}\frac{T_{2n}(v_p)-1}{v_p-1},$
\end{enumerate}
where $u_p$ and $v_p$ are the same as in the statement of the theorem. Hence,
$$\tau(n)=(-1)^{(s_k-1)}\frac{n}{2}\prod_{p=1}^{s_k-1}\frac{T_n(u_p)-1}{u_p-1}
\prod_{p=1}^{s_k}\frac{T_{2n}(v_p)-1}{T_n(v_p)-1}$$
We note that $P^{\prime}(1)=-2(s_1^2+s_2^2+\ldots+s_k^2).$ Since $P(1)=1$ and $P(u)$
is the order $s_k$ polynomial whose leading term is $-2^{s_k}u^{s_k},$ we have
$\prod_{p=1}^{s_k-1}(u_p-1)=(-2)^{-s_k}P^{\prime}(1),$ where $u_1,u_2,\ldots,u_{s_k-1}$
are all roots of the equation $P(u)-1=0$ different from $1.$ Finally, taking into account
these properties and the identity $T_{2n}(v_p)-1=2(T_n(v_p)-1)(T_n(v_p)+1)$ we obtain
$$\tau(n)=\frac{n\,4^{s_k-1}}{q}\prod_{p=1}^{s_k-1}(T_n(u_p)-1)\prod_{p=1}^{s_k}(T_n(v_p)+1),$$
where $q=s_1^2+s_2^2+\ldots+s_k^2.$

\section{Arithmetic properties of the complexity for circulant graphs}\label{circarithm}

It was noted in the series of paper (\cite{XiebinLinZhang}, \cite{ZhangYongGol}, 
\cite{ZhangYongGolin}) that in many   cases the complexity of circulant graphs
is given by the formula $\tau(n)=n a(n)^2,$ where $a(n)$ is an integer sequence. In the
same time, this is not always true. Indeed, for the graph $C_n(1,3)$ and $n$ even we have
$\tau(n)=2 n a(n)^2$ for some integer sequence $a(n).$

The aim of the next theorem is to explain this phenomena. Recall that any positive integer
$p$ can be uniquely represented in the form $p=q \,r^2,$ where $p$ and $q$ are positive
integers and $q$ is square-free. We will call $q$ the \textit{square-free part} of $p.$

\begin{theorem}\label{lorenzini}
Let $\tau(n)$ be the number of spanning trees in the circulant graph
$C_{n}(s_1,s_2,\ldots,s_k),${\ }$\,1\le s_1<s_2<\ldots<s_k<\frac{n}{2}.$
Denote by $p$ the number of odd elements in the sequence $s_1,s_2,\ldots,s_k$
and let $q$ be the square-free part of $p.$
Then there exists an integer sequence $a(n)$ such that
\begin{enumerate}
\item[ $1^0$ ]  $\tau(n)=n\,a(n)^2,$ if   $n$ is odd;
\item[ $2^0$ ]  $\tau(n)=q \,n\,a(n)^2,$  if $n$ is even.
\end{enumerate}
\end{theorem}
\begin{pf} The number of odd elements in the sequence $s_1,s_2,s_3,\ldots,s_k$ is counted by the formula
$p=\sum\limits_{i=1}^{k}\frac{1-(-1)^{s_i}}{2}.$ If $n$ is even and the graph $C_{n}(s_1,s_2,s_3,\ldots,s_k )$
is connected then at least one of the numbers $s_1,s_2,s_3,\ldots,s_k$ is odd, otherwise the number of spanning
trees $\tau(n)=0.$ So, we can assume that $p>0.$

We already know that all non-zero eigenvalues of the graph $C_{n}(s_1,s_2,s_3,\ldots,s_k )$
are given by the formulas
$\lambda_j=L(\varepsilon_n^{j}),\,j=1,\ldots,n-1,$ where $L(z)=2k-\sum\limits_{i=1}^{k}(z^{s_i}+z^{-s_i})$ and
$\varepsilon_n=e^{\frac{2\pi i}{n}}.$ We note that $\lambda_{n-j}=L(\varepsilon_n^{n-j})=L(\varepsilon_n^{j})=\lambda_j.$

By the Kirchhoff theorem we have $n\,\tau(n)=\prod\limits_{j=1}^{n-1}\lambda_j.$
Since $\lambda_{n-j}=\lambda_j,$  we obtain
$n\,\tau(n)=(\prod\limits_{j=1}^{\frac{n-1}{2}}\lambda_j)^2$ if $n$ is odd and
$n\,\tau(n)=\lambda_{\frac{n}{2}}(\prod\limits_{j=1}^{\frac{n}{2}-1}\lambda_j)^2$
if $n$ is even. We note that each algebraic number $\lambda_j$ comes with all its
Galois conjugate \cite{Lor}. So, the numbers
$c(n)=\prod\limits_{j=1}^{\frac{n-1}{2}}\lambda_j$ and $d(n)=\prod\limits_{j=1}^{\frac{n}{2}-1}\lambda_j$
are integers. Also, for even $n$ we have $\lambda_{\frac{n}{2}}=2k-\sum\limits_{i=1}^{k}
((-1)^{s_i}+(-1)^{-s_i})=2\sum\limits_{i=1}^{k}(1-(-1)^{s_i})=4p.$
Hence, $n\,\tau(n)= c(n)^2$ if $n$ is odd and $n\,\tau(n)=4 p\,d(n)^2$ if $n$ is even.
Let $q$ be the free square part of $p$ and $p=q\,r^2.$ The circulant graph
$C_{n}(s_1,s_2,s_3,\ldots,s_k )$ has a cyclic group of automorphisms $\mathbb{Z}_n$
acting fixed point free on the set of all spanning trees, therefore $\tau(n)$ is a
multiple of $n.$  As a result, for an integer number $\frac{\tau(n)}{n}$ we have
\begin{enumerate}
\item $\displaystyle{\frac{\tau(n)}{n}=\left(\frac{c(n)}{n}\right)^2}$  if $n$  is odd and
\item  $\displaystyle{\frac{\tau(n)}{n}=q \left(\frac{2\,r\,d(n)}{n}\right)^2}$ if $n$  is even.
\end{enumerate} Setting $a(n)=\frac{c(n)}{n}$ in the first case and $a(n)=\frac{2\,r\,d(n)}{n}$
in the second, we conclude that number $a(n)$ is always integer and the statement of theorem follows.
\end{pf}
\bigskip

The following theorem clarifies some number-theoretical properties of the complexity $\tau(n)$
for circulant graphs of odd valency.
\bigskip

\begin{theorem}\label{lorenzinin} Let $\tau(n)$ be the number of spanning trees in the circulant graph
$$C_{2n}(s_1,s_2,s_3,\ldots,s_k,n),\,1\le s_1< s_2<\ldots< s_k<n.$$ Denote by $p$ the number
of odd elements in the sequence $s_1,s_2,s_3,\ldots,s_k.$ Let $q$ be the square-free part of $2p$ and
$r$ be the square-free part of $2p+1.$ Then there exists an integer sequence $a(n)$ such that
\begin{enumerate}
\item[$1^0.$]  $\tau(n)=r\,n\,a(n)^2,$ if $n$ is odd;
\item[$2^0.$]  $\tau(n)=q\,n\,a(n)^2,$ if $n$ is even.
\end{enumerate}
\end{theorem}
\begin{pf} The number $p$ of odd elements in the sequence $s_1,s_2,\ldots,s_k$ is counted by the formula
$p=\sum\limits_{i=1}^{k}\frac{1-(-1)^{s_i}}{2}.$ The non-zero eigenvalues of the graph $C_{2n}(s_1,s_2,\ldots,s_k,n)$
are given by the formulas $\lambda_j=L(\varepsilon_{2n}^{j})-(-1)^j,\,j=1,2,\ldots,2n-1,$ where
$L(z)=2k+1-\sum\limits_{l=1}^{k}(z^{s_l}+z^{-s_l})$ and $\varepsilon_{2n}=e^{\frac{\pi i}{n}}.$

By the Kirchhoff theorem we have $2n\,\tau(n)=\prod\limits_{j=1}^{2n-1}\lambda_j.$ Since $\lambda_{2n-j}=\lambda_j,$ we
obtain $2n\,\tau(n)=\lambda_{n}(\prod\limits_{j=1}^{n-1}\lambda_j)^2,$ where $\lambda_n=L(-1)-(-1)^n.$ Now  we have
$$\lambda_n=2k+1-(-1)^n-2\sum\limits_{l=1}^{k}(-1)^{s_l}=1-(-1)^n+4\sum\limits_{l=1}^{k}\frac{1-(-1)^{s_l}}{2}=1-(-1)^n+4p.$$
So, $\lambda_n=4\,p,$ if $n$ is even and  $\lambda_n=4\,p+2,$ if $n$ is odd. We note that each algebraic number $\lambda_j$
comes in $\prod\limits_{j=1}^{n-1}\lambda_j$ together with all its Galois conjugate, so the number
$c(n)=\prod\limits_{j=1}^{n-1}\lambda_j$ is an integer \cite{Lor}.

Hence, $n\,\tau(n)=(2\,p+1)c(n)^2,$ if $n$ is odd and $n\,\tau(n)=2 p\,c(n)^2,$ if $n$ is even.
Let $q$ and $r$ be the free square parts of $2p$ and of $2p+1$ respectively. Then for some
integers $x$ and $y$ we have $2p=q\,x^2$ and $2\,p+1=r\,y^2.$ The circulant graph
$C_{2n}(s_1,s_2,s_3,\ldots,s_k,n)$ has a cyclic group of automorphisms $\mathbb{Z}_n$ acting
fixed point free on the set of all spanning trees, therefore $\tau(n)$ is a multiple of $n.$

It is important to note that the cyclic group of automorphisms $\mathbb{Z}_{2n}$ acts not
fixed point free on the set of all spanning trees of $C_{2n}(s_1,s_2,\ldots,s_k, n ).$ Some
trees whose edges joint apposite vertices of the graph are fixed by the involution from
$\mathbb{Z}_{2n}.$ So, the $\tau(n)$ is not necessary divided by $2n.$

Now, the integer number $\frac{\tau(n)}{n}$ can be represented in the form
\begin{enumerate}
\item $\displaystyle{\frac{\tau(n)}{n}=r\,\left(\frac{x\,c(n)}{n}\right)^2}$  if $n$  is odd and
\item $\displaystyle{\frac{\tau(n)}{n}=q \left(\frac{y\,c(n)}{n}\right)^2}$ if $n$  is even.
\end{enumerate} Setting $a(n)=\frac{x\,c(n)}{n}$ in the first case and $a(n)=\frac{y\,c(n)}{n}$
in the second, we conclude that  number $a(n)$ is always integer. The theorem is proved.
\end{pf}

\section{Asymptotic for the number of spanning trees}\label{assection}

In this section we give asymptotic formulas for the number of spanning trees in circulant
graphs. It is interesting to compare these results with those from papers \cite{Louis}, \cite{XiebinLinZhang}, \cite{ZhangYongGol},
 and  \cite{ZhangYongGolin}, where the similar results were
obtained by different methods.

\begin{theorem}\label{asymptotic1} Let ${\rm gcd}\,(s_1,s_2,\ldots,s_k)=d.$ Then
the number of spanning trees in the circulant graph
$C_{n}(s_1,s_2,\ldots,s_k),\,1\le s_1< s_2<\ldots< s_k<\frac{n}{2}$
has the following asymptotic
$$ \tau(n)\sim\frac{n\,d^2 }{q}A^n,\text{ as }n\to\infty \text{ and }(n,d)=1,$$
where $q=s_1^2+s_2^2+\ldots+s_k^2$ and $A=\exp(\int_0^1\log|L(e^{2\pi i t})|dt)$
is the Mahler measure of Laurent polynomial $L(z)=2k -\sum\limits_{i=1}^k(z^{s_i}+z^{-{s_i}}).$

\end{theorem}
\begin{pf} Since we are interested only in connected graphs, the condition $(n,d)=1$ is always
satisfied.  Then the graphs $C_{n}(s_1,s_2,\ldots,s_k)$ and $C_{n}(s_1/d,s_2/d,\ldots,s_k/d)$
are isomorphic and it is sufficient to prove the theorem for the case $d=1.$ By Corollary \ref{corollary21}, 
the number of spanning trees $\tau(n)$ is given by
$$\tau(n)=\frac{n}{q}\prod_{p=1}^{s_k-1}|2\,T_n(w_p)-2|.$$

By Lemma~\ref{lemma3}, we have $T_n(w_s)=\frac{1}{2}(z_s^n+z_s^{-n}),$ where the $z_s$ and $1/z_s$
are roots of the polynomial $L(z)$ with the property $|z_s|\neq1,\,s=1,2,\ldots,s_k-1.$ Replacing
$z_s$ by $1/z_s,$ if it is necessary, we can assume that $|z_s|>1$ for all $s=1,2,\ldots,s_k-1.$
Then $T_n(w_s)\sim\frac{1}{2}z_s^n,$ as $n$ tends to $\infty.$ So, 
$|2T_n(w_s)-2|\sim|z_s|^n,\,\,n\to\infty.$
Hence
$$\prod_{s=1}^{s_k-1}|2\,T_n(w_s)-2|\sim \prod_{s=1}^{s_k-1}|z_s|^n=
\prod\limits_{L(z)=0,\,|z|>1}|z|^n=A^n,$$
where $A=\prod\limits_{L(z)=0,\,|z|>1}|z|$  is the Mahler measure of $L(z).$
By the results mentioned in the preliminary part, it can be found by the formula
$A=\exp(\int_0^1\log|L(e^{2\pi i t})|dt).$

Finally, $$\tau(n)=\frac{n}{q}\prod_{s=1}^{s_k-1}|2\,T_n(w_s)-2|\sim\frac{n}{q}A^n,\,n\to \infty.$$

If $d>1,$ one has to substitute the sequence $s_1,s_2,\ldots,s_k$ by $s_1/d,s_2/d,\ldots,s_k/d.$
In this case, the associated polynomial is $L_d(z)=2k-\sum\limits_{i=1}^k(z^{s_i/d}+z^{-s_i/d}).$
Since $L_d(z^d)=L(z)$ and $M(L_d(z))=M(L(z)),$ the statement of the theorem follows.
\end{pf}
\bigskip

As an immediate consequence of Theorem~\ref{asymptotic1} we have the following result.
\bigskip

\begin{corollary}\label{cor1}
The thermodynamic limit of the graph sequence $C_{n}(s_1,s_2,\ldots,s_k)$  is equal to the small
Mahler measure of Laurent polynomial $L(z)=2k-\sum\limits_{i=1}^k(z^{s_i}+z^{-s_i}).$   That is
$$\lim\limits_{n\to\infty}\frac{\log\tau(C_{n}(s_1,s_2,\ldots,s_k,n))}{n}=m(L),$$
where $m(L)=\int\limits_0^1\log|L(e^{2\pi i t})|dt.$
\end{corollary}

\bigskip

The next theorem is a direct consequence of Theorem~\ref{odddegree} and can be proved by  the same arguments as Theorem~\ref{asymptotic1}.


\bigskip

\begin{theorem}\label{asymptotic3} Let ${\rm gcd}\,(s_1,s_2,\ldots,s_k)=d.$ Then the number of spanning trees in the circulant graph
$C_{2n}(s_1,s_2,\ldots,s_k,n),\,1\le s_1< s_2<\ldots< s_k<n$ has the following asymptotic
$$\tau(n)\sim\frac{ n\,d^2 }{2\,q}K^n,\text{ as }n\to\infty \text{   and }  (n,d)=1.$$
Here $q=s_1^2+s_2^2+\ldots+s_k^2$ and $K=\exp(\int\limits_0^1\log|L^2(e^{2\pi i t})+2 L(e^{2\pi i t})|dt)$
is the Mahler measure of the Laurent polynomial $L(z)(L(z)+2),$ where $L(z)=2k-\sum\limits_{i=1}^k(z^{s_i}+z^{-s_i}).$

\end{theorem}
\bigskip

As a   corollary  of Theorem~\ref{asymptotic3}  we have the following result.

\bigskip
\begin{corollary}\label{cor2} The thermodynamic limit of the sequence $C_{2n}(s_1,s_2,\ldots,s_k,n)$ of circulant graphs
is equal to the small Mahler measure of Laurent polynomial $R(z)=L(z)(L(z)+2),$ where $L(z)=2k-\sum\limits_{i=1}^k(z^{s_i}+z^{-s_i}).$ More precisely,
$$\lim\limits_{n\to\infty}\frac{\log\tau(C_{2n}(s_1,s_2,\ldots,s_k,n))}{n}=m(R),$$
where $m(R)=\int\limits_0^1\log|R(e^{2\pi i t})|dt.$
\end{corollary}

\section{Examples}\label{tables}

\begin{enumerate}
\item{\textbf{Graph $C_n(1,2)$}.}
From the paper \cite{BoePro} we have $\tau(n)=nF_n^2,$ where $F_n$ is the $n$-th Fibonacci number. By Theorem~\ref{asymptotic1}, for $A_{1,2}=\frac {1}{2}(3+\sqrt{5})$
and $q=5$ we have the following asymptotic $\tau(n)\sim\frac{n}{5}A_{1,2}^n,\, n\to\infty.$

\item{\textbf{Graph $C_n(1,3)$}.}
By Theorem~\ref{lorenzini}, we have $\tau(n)=n\,a(n)^2$ in $n$  is odd, and
$\tau(n)=2n\,a(n)^2$ in $n$ is even, where $a(n)$ is an integer sequence. As a consequence of Theorem~3 in \cite{YTA97},
one can show that $a(n)$ is
\textit{A112835} sequence in the On - Line Encyclopedia of Integer Sequences.

In this case, $A_{1,3}=\frac12(1+\sqrt{1-2i})(1+\sqrt{1+2i})\approx2.89,\,q=10$ and
$\tau(n)\sim\frac{n}{10}A_{1,3}^n,\, n\to\infty.$  (Compare with Example 2 in \cite{Xiebin}.)

\item{\textbf{Graph  $C_n(2,3)$}.}
By Theorem~\ref{lorenzini}, we have $\tau(n)=n\,a(n)^2$ for some integer sequence $a(n).$ One can check
(see, for example \cite{ZhangYongGol}, Theorem~9) that $a(n)$ satisfies the linear recursive relation
$a(n)= a(n-1)+a(n-2)+a(n-3)-a(n-4)$ with initial data $a(0)=0, \,a(1)=1, \,a(2)=1, \,a(3)=1.$ Note $a(n)$ is
\textit{A116201} sequence in the On - Line Encyclopedia of Integer Sequences.

In this case $A_{2,3}\approx2.96$ is the root of the equation $1 - 3 z + z^2 - 3 z^3 + z^4=0, \,q=13$  and
$\tau(n)\sim\frac{n}{13}A_{2,3}^n,\, n\to\infty.$ (See also \cite{Xiebin}, Example 3.)

\item{\textbf{Graph  $C_n(1,2,3)$}.}
Here $A_{1,2,3}=\frac12(2+\sqrt{7}+\sqrt{7+4\sqrt{7}})\approx4.42$ and $\,\tau(n)\sim\frac{n}{14}A_{1,2,3}^n,\, n\to\infty.$
By Theorem~\ref{lorenzini}, there exists an integer sequence $a(n)$ such that $\tau(n)= n\,a(n)^2$ if $n$ is odd,
and $\tau(n)=2n\,a(n)^2$ if $n$ is even. (Compare with \cite{Xiebin}, Example 4.)

\item{\textbf{Graph M\"obius ladder $C_{2n}(1,n).$}}
In this case, by \cite{BoePro} we have $\tau(n)=n(T_n(2)+1)\sim\frac{n}{2}(2+\sqrt{3})^n,\,n\to\infty.$
Also, by Theorem~\ref{lorenzinin}, there exists an integer sequence $a(n)$ such that $\tau(n)=3 n\,a(n)^2$ if $n$ is odd,
and $\tau(n)=2n\,a(n)^2$ if $n$ is even. By Corollary~4   from \cite{MedMed2}  one can conclude that
$a(2m+1)=T_m(2)+U_{m-1}(2)$ and $a(2m)=T_m(2).$

\item{\textbf{Graph $C_{2n}(2,n),\,n$ is odd (the $n$-prism graph).}}
The number of spanning trees  $\tau(n)=n(T_n(2)-1)\sim\frac{n}{2}(2+\sqrt{3})^n,\,n\to\infty.$
Also, by Theorem~\ref{lorenzinin}, there is an integer sequence $a(n)$ such that $\tau(n)= n\,a(n)^2$ if $n$ is odd.
Indeed, for $n=2m+1$ we have $a(2m+1)=T_m(2)+3 U_{m-1}(2).$

\item{\textbf{Graph  $C_{2n}(1,2,n)$}.}

$K_{1,2}=\frac14(3+\sqrt{5})(4+\sqrt{3}+\sqrt{15+8\sqrt{3}})\approx14.54,\,\tau(n)\sim\frac{n}{10}\,K_{1,2}^n,\, n\to\infty.$

By Theorem~\ref{lorenzinin}, there exists an integer sequence $a(n)$ such that $\tau(n)=3n\,a(n)^2$ if $n$ is odd and
$\tau(n)=2n\,a(n)^2$ if $n$ is even.

\item{\textbf{Graph  $C_{2n}(1,2,3,n)$}.}
$K_{1,2,3}\approx 32.7865,\,\tau(n)\sim\frac{n}{28}K_{1,2,3}^{n},\,n\to\infty.$

By Theorem~\ref{lorenzinin}, for some integer sequence $a(n)$ we have $\tau(n)=5n\,a(n)^2$ if $n$ is odd and
$\tau(n)=n\,a(n)^2$ if $n$ is even.
\end{enumerate}

\section*{ACKNOWLEDGMENTS}
This work was supported by by the Russian Foundation for Basic Research (projects 15-01-07906, and 16-31-00138)
and the Slovenian-Russian grant (2016--2017).

\end{document}